\documentclass[12pt, a4paper]{article}

\usepackage{amsfonts}
\usepackage{amssymb}
\usepackage{amscd}
\usepackage{mathrsfs}
\usepackage{amsmath,amssymb, amsthm}
\usepackage{graphicx}
\usepackage{color}
\usepackage{authblk}
\usepackage{indentfirst,latexsym,amssymb}
\usepackage[bookmarks=false]{hyperref}

\theoremstyle{definition}

\newtheorem*{exam}{Example}

\def\Cay{\hbox{\rm Cay}}

\usepackage{enumerate}
\long\def\delete#1{}

\usepackage{xcolor}
\usepackage[normalem]{ulem}

\input amssym.def
\input amssym.tex

\newcommand{\be}{\begin{equation}}
\newcommand{\ee}{\end{equation}}
\newcommand{\bea}{\begin{eqnarray}}
\newcommand{\eea}{\end{eqnarray}}
\newcommand{\bean}{\begin{eqnarray*}}
\newcommand{\eean}{\end{eqnarray*}}

\newtheorem{manualtheoreminner}{Theorem}
\newenvironment{manualtheorem}[1]{%
  \renewcommand\themanualtheoreminner{#1}%
  \manualtheoreminner
}{\endmanualtheoreminner}

\newtheorem{manuallemmainner}{Lemma}
\newenvironment{manuallemma}[1]{%
  \renewcommand\themanuallemmainner{#1}%
  \manuallemmainner
}{\endmanuallemmainner}

\newtheorem{manualcoroinner}{Corollary}
\newenvironment{manualcoro}[1]{%
  \renewcommand\themanualcoroinner{#1}%
  \manualcoroinner
}{\endmanualcoroinner}

\title{Corrigendum to ``On subgroup perfect codes in Cayley graphs'' [European J. Combin. 91 (2021)
103228]}

\author[a]{Junyang Zhang}
\author[b]{Sanming Zhou}
\affil[a]{{\small School of Mathematical Sciences, Chongqing Normal University, Chongqing 401331, P. R. China}}
\affil[b]{{\small School of Mathematics and Statistics, The University of Melbourne, Parkville, VIC 3010, Australia}}

\date{}

\begin{document}

\openup 0.5\jot
\maketitle

\renewcommand{\thefootnote}{\empty}%{footnote}}
\footnotetext{E-mail addresses: jyzhang@cqnu.edu.cn (Junyang Zhang), sanming@unimelb.edu.au (Sanming Zhou)}

\begin{abstract}
We correct the statements of two theorems and two corollaries in our paper [On subgroup perfect codes in Cayley graphs, \emph{European J. Combin.} 91 (2021) 103228]. Proofs of these theorems and three other results are given as well.
\end{abstract}

\delete
{
}


This is a corrigendum to our paper [On subgroup perfect codes in Cayley graphs, \emph{European J. Combin.} 91 (2021) 103228], which will be referred to as \cite{ZZ21} hereafter. After the paper was published, an error was found in the proof of Theorem \ref{basic}, thereby rendering the result invalid. Here we will give a corrected version of this result and its proof. As a consequence of this change, we need to revise \cite[Theorem 3.2]{ZZ21} and its proof as well. It turns out that the subgroup involved in Theorem 3.2 and Corollaries 3.3 and 3.4 in \cite{ZZ21} is required to be a $2$-group or have odd order or odd index. The proofs of Theorems \ref{quotient}, \ref{ns} and \ref{QK} in \cite{ZZ21} need to be modified due to their reliance on Corollaries \ref{equivalent} and \ref{equivalent2}, but the statements of these results remain true. We will give new proofs of these three theorems using the corrected Theorem 3.1.

We will use the notation and terminology in \cite{ZZ21}. As usual, denote by $|G|$ the order of a group $G$, $|G:H|$ the index in $G$ of a subgroup $H$ of $G$, and $|X|$ the cardinality of a set $X$. We will also use the same labels to label corresponding results from \cite{ZZ21} or their revised versions.

The following four lemmas are exactly Lemmas 2.1--2.4 in \cite{ZZ21}. We repeat them here for later convenience.

\begin{manuallemma}{2.1}
\label{Ct}
{\em Let $G$ be a group and $H$ a subgroup of $G$. Then $H$ is a perfect code of $G$ if and only if it has a Cayley transversal in $G$.}
\end{manuallemma}

\begin{manuallemma}{2.2}
\label{sub}
{\em Let $G$ be a group and $H$ a subgroup of $G$. Then $H$ is a perfect code of $G$ if and only if it is a perfect code of any subgroup of $G$ which contains $H$.}
\end{manuallemma}

\begin{manuallemma}{2.3}
\label{conjugate}
{\em Let $G$ be a group and $H$ a subgroup of $G$. If $H$ is a perfect code of $G$, then for any $g \in G$, $g^{-1}Hg$ is a perfect code of $G$. More specifically, if $H$ is a perfect code in $\Cay(G, S)$ for some connection set $S$ of $G$, then $g^{-1}Hg$ is a perfect code in $\Cay(G, g^{-1}Sg)$.}
\end{manuallemma}

\begin{manuallemma}{2.4}
\label{HXZ}
{\em (\cite[Theorem 2.2]{HXZ18}) Let $G$ be a group and $H$ a normal subgroup of $G$. Then $H$ is a perfect code of $G$ if and only if for all $x\in G$, $x^{2}\in H$ implies $(xh)^{2}=1$ for some $h\in H$.}
\end{manuallemma}

An error occurs in lines 11--12 of the proof of \cite[Theorem 3.1]{ZZ21}, where the claim that ``$xH\neq yH$ for distinct elements $x, y \in X$'' is not true as $z^{-1}H$ may be equal to $wH$ for some $w\in T$. The corrected version of \cite[Theorem 3.1]{ZZ21} should read as follows.

\begin{manualtheorem}{3.1}
\label{basic}
{\em Let $G$ be a group and $H$ a subgroup of $G$. Then $H$ is not a perfect code of $G$ if and only if
there exists a double coset $D=HxH$ with $D=D^{-1}$ having an odd number of left cosets of $H$ in $G$ and containing no involution.
In particular, if $H$ is not a perfect code of $G$, then there exists a $2$-element $x\in G\setminus H$ such that $x^{2}\in H$, $|H:H\cap xHx^{-1}|$ is odd, and $HxH$ contains no involution.}
\end{manualtheorem}

\begin{proof}
We prove the sufficiency first. Assume that there is a double coset $D=HxH$ with $D=D^{-1}$ having $m$ left cosets of $H$ in $G$ for some odd integer $m$ and containing no involution. Set $D=\cup_{i=1}^{m}g_{i}H$, where $g_{1}H, \ldots, g_{m}H$ are distinct left cosets of $H$ in $G$. Suppose to the contrary that $H$ is a perfect code of $G$. By Lemma \ref{Ct}, $H$ has a Cayley transversal $T$ in $G$. It follows that for each $1\leq i\leq m$, there exists $x_i \in G$ such that $T\cap g_{i}H = \{x_{i}\}$. Set $X=\{x_1,x_2,\ldots,x_m\}$. Then $X=T\cap D$. Since both $T$ and $D$ are inverse-closed, we have
$X^{-1}=(T\cap D)^{-1}=T^{-1}\cap D^{-1}=T\cap D=X$. Since $m$ is odd, $X$ contains at least one involution. Since $X$ is a subset of $D$, it follows that $D$ contains at least one involution, but this contradicts our assumption.

Now we prove the necessity. Assume that $H$ is not a perfect code of $G$.
Take a subset $T$ of $G$ with maximum cardinality such that $1\in T$, $T^{-1}=T$, $HTH=TH$ and $xH\neq yH$ for all pairs of distinct elements $x,y\in T$. (The existence of $T$ follows from the fact that there are subsets of $G$, say, $\{1\}$, with all these properties.) Since $H$ is not a perfect code of $G$,
by Lemma \ref{Ct}, $T$ is not a left transversal of $H$ in $G$. It follows that  $G \setminus TH \neq\emptyset$ and therefore we can take an element $x\in G \setminus TH$. Set $D=HxH$. Then $D=Hx^{-1}H$.

Since $HTH=TH$ and $x\notin TH$, we have $D\cap TH=\emptyset$. Furthermore, since $T=T^{-1}$, we have
\begin{equation*}
D^{-1}\cap TH=D^{-1}\cap HTH=(D\cap HTH)^{-1}=\emptyset.
\end{equation*}
Set $|H:H\cap xHx^{-1}|=\ell$. Since $H\cap x^{-1}Hx=x^{-1}(H\cap xHx^{-1})x$, we have $|H:H\cap x^{-1}Hx|=\ell$. It is straightforward to verify that $h_{1}xH=h_{2}xH$ if and only if $h_{1}(H\cap xHx^{-1})=h_{2}(H\cap xHx^{-1})$ for any pair of elements $h_{1},h_{2}\in H$. Therefore, $D$ is the union of $\ell$ distinct left cosets of $H$ in $G$. So we can express $D$ as the union of $\ell$ distinct left cosets $x_{1}xH,\ldots,x_{\ell}xH$ for some $x_{1},\ldots,x_{\ell} \in H$. Similarly, we can write $D^{-1}$ as the union of $\ell$ distinct left cosets $y_{1}x^{-1}H,\ldots,y_{\ell}x^{-1}H$ for some $y_{1},\ldots,y_{\ell}\in H$.

We will complete the proof of the necessity by showing that $D=D^{-1}$, $D$ contains no involution, and $\ell$ is odd.
If $D\neq D^{-1}$, then we set
\begin{equation*}
X=T\cup\{x_{1}xy_{1}^{-1},\ldots,x_{\ell}xy_{\ell}^{-1},y_{1}x^{-1}x_{1}^{-1},
\ldots,y_{\ell}x^{-1}x_{\ell}^{-1}\}.
\end{equation*}
If $D$ contains an involution $z$, then $z\in x_{i}xH$ for some
$1\leq i\leq\ell$. We set
\begin{equation*}
X=T\cup\{x_{1}x_{i}^{-1}zx_{i}x_{1}^{-1},\ldots,
x_{\ell}x_{i}^{-1}zx_{i}x_{\ell}^{-1}\}.
\end{equation*}
If $D=D^{-1}$ and $\ell$ is even, then $x^{-1}H=hxH$ for some $h\in H$. We set
\begin{equation*}
X=T\cup \left\{x_{1}xhx_{\frac{\ell}{2}+1}^{-1},\ldots,x_{\frac{\ell}{2}}xhx_{\ell}^{-1},
x_{\frac{\ell}{2}+1}h^{-1}x^{-1}x_{1}^{-1},\ldots,x_{\ell}h^{-1}x^{-1}x_{\frac{\ell}{2}}^{-1}\right\}.
\end{equation*}
In each case above, we have defined a subset $X$ of $G$ which contains $T$ as a proper subset. It can be verified that $X^{-1} = X$, $HXH=XH$ and $xH\neq yH$ for any pair of distinct elements $x,y\in X$, but this contradicts the maximality of $T$. This contradiction shows that $D=D^{-1}$, $D$ contains no involution, and $\ell$ is an odd integer.

It remains to prove the last statement in the theorem. Suppose that $H$ is not a perfect code of $G$. Then by what we have proved above there exists $g\in G\setminus H$ such that $Hg^{-1}H=HgH$, $|H:H\cap gHg^{-1}|$ is odd, and $HgH$ contains no involution. Since $Hg^{-1}H=HgH$, we have $Hg^{-1}=Hgh$ for some $h\in H$. Set $y=gh$. Then $y\in G\setminus H$ and $y^2\in H$. Let $s$ be the largest odd divisor of the order of $y$. Set $x=y^{s}$. Then $x$ is a $2$-element and $x^2 \in H$. Since $s-1$ is even and $y^2\in H$, we have $y^{s-1}\in H$. Hence $HxH=HyH=HgH$ and $xHx^{-1}=yHy^{-1}=ghHh^{-1}g^{-1}=gHg^{-1}$. Note that $x \not \in H$ as $y \not \in H$. Therefore, $x \in G\setminus H$ is a $2$-element such that $x^{2}\in H$,  $|H:H\cap xHx^{-1}|$ is odd, and $HxH$ contains no involution.
\end{proof}

The corrected version of \cite[Theorem 3.2]{ZZ21} should read as follows.

\begin{manualtheorem}{3.2}
\label{normal}
{\em Let $G$ be a group and $H$ a subgroup of $G$. Suppose that either $H$ is a $2$-group or at least one of $|H|$ and $|G:H|$ is odd. Then $H$ is a perfect code of $G$ if and only if $H$ is a perfect code of $N_{G}(H)$.}
\end{manualtheorem}

\begin{proof}
The necessity follows from Lemma \ref{sub}. To prove the sufficiency, we assume that $H$ is a perfect code of $N_{G}(H)$. By way of contradiction, suppose that $H$ is not a perfect code of $G$. By Theorem \ref{basic}, there exists a $2$-element $x\in G\setminus H$ such that $x^2\in H$, $|H:H\cap xHx^{-1}|$ is odd, and $HxH$ contains no involution. In particular, $xH$ contains no involution. Since $H$ is a perfect code of $N_{G}(H)$, we have $x\notin N_{G}(H)$ by Lemma \ref{HXZ}. Set $L=H\cap xHx^{-1}$. Then $L$ is a proper subgroup of $H$. Since $|H:L|$ is odd, $H$ cannot be a $2$-group. Since $HxH$ contains no involution, $x$ is not an involution. Therefore, $x^2$ is a non-identity $2$-element. Since $x^2\in H$, it follows that $H$ is of even order. Since
\begin{equation*}
xLx^{-1} = x(H\cap xHx^{-1})x^{-1} = xHx^{-1}\cap x^2Hx^{-2} = xHx^{-1}\cap H = L,
\end{equation*}
we have $x\in N_{G}(L)$. It follows that $|N_{G}(L):L|$ is even. Since $|N_{G}(L):L|$ is a divisor of $|G:L|$, it follows that $|G:L|$ is even. Since $|G:L|=|G:H||H:L|$ and $|H:L|$ is odd, we conclude that $|G:H|$ is even. Now we have proved that $H$ is not a $2$-group and both $|H|$ and $|G:H|$ are even. This contradicts our assumption and therefore $H$ must be a perfect code of $G$.
\end{proof}

The following is a counterexample to both \cite[Theorem 3.1]{ZZ21} and \cite[Theorem 3.2]{ZZ21}. 

\begin{exam}
\label{example}
The special linear group of degree $2$ over the finite field of order $3$ has the following representation
\begin{equation*}
G=\langle a,b,c\mid a^4=c^3=1, a^2=b^2, ab=b^{-1}a, ac=cb, bc=cab\rangle.
\end{equation*}
Note that $G$ is a split extension of its Sylow $2$-subgroup $\langle a,b\rangle$ by its Sylow $3$-subgroup $\langle c\rangle$. Note also that $\langle a,b\rangle$ is isomorphic to the quaternion group, of which $a^2$ is the unique involution. Therefore $a^2$ is the unique involution of $G$. Set $H=\langle a^2,c\rangle$. Then $H$ is a cyclic subgroup of $G$ of order $6$. Clearly, $a\in G\setminus H$ and $a^{2}\in H$. Furthermore, it is straightforward to verify that $HaH=Ha^{-1}H=G\setminus H$. Therefore, $HaH$ is an inverse-closed double coset having $3$ left cosets of $H$ in $G$ and containing no involution. By Theorem \ref{basic}, $H$ is not a perfect code of of $G$. Since $aH=a^{-1}H$ and $aH$ contains no involution, this shows that the statement in \cite[Theorem 3.1]{ZZ21} is not true. Moreover, we have $H=N_{G}(H)$ and hence $H$ is a perfect code of $N_{G}(H)$. This shows that \cite[Theorem 3.2]{ZZ21} is not true without the additional condition that either $H$ is a $2$-group or at least one of $|H|$ and $|G:H|$ is odd.
\end{exam}

Combining Lemma \ref{HXZ} and Theorem \ref{normal}, we obtain the following corrected version of \cite[Corollary 3.3]{ZZ21}.

\begin{manualcoro}{3.3}
\label{equivalent}
{\em Let $G$ be a group and $H$ a subgroup of $G$. Suppose that either $H$ is a $2$-group or at least one of $|H|$ and $|G:H|$ is odd. Then $H$ is a perfect code of $G$ if and only if for any $x\in N_{G}(H)$, $x^{2}\in H$ implies $(xh)^{2}=1$ for some $h\in H$.}
\end{manualcoro}

The corrected version of \cite[Corollary 3.4]{ZZ21} should read as follows. The proof of this new version is exactly the same as the proof of \cite[Corollary 3.4]{ZZ21}, except that Corollary 3.3 should be understood as the corrected Corollary \ref{equivalent}.

\begin{manualcoro}{3.4}
\label{equivalent2}
{\em Let $G$ be a group and $H$ a subgroup of $G$. Suppose that either $H$ is a $2$-group or at least one of $|H|$ and $|G:H|$ is odd. Then $H$ is a perfect code of $G$ if and only if for any $2$-element $x\in N_{G}(H)$, $x^{2}\in H$ implies $(xh)^{2}=1$ for some $h\in H$.}
\end{manualcoro}

The proof of \cite[Theorem 3.7]{ZZ21} relies on \cite[Corollary 3.3]{ZZ21}. Since the latter has been modified, we have to revise the proof of the former. Below we will give a proof of \cite[Theorem 3.7]{ZZ21} using Lemmas \ref{Ct} and \ref{HXZ} and the corrected Theorem \ref{basic}. Fortunately, both statements in \cite[Theorem 3.7]{ZZ21} remain true under the same hypothesis. We repeat this theorem here for convenience.

\begin{manualtheorem}{3.7}
\label{quotient}
{\em Let $G$ be a group, $N$ a normal subgroup of $G$, and $H$ a subgroup of $G$ which contains $N$. Then the following hold:
\begin{enumerate}[\rm (a)]
  \item if $H$ is a perfect code of $G$, then $H/N$ is a perfect code of $G/N$;
  \item if $N$ and $H/N$ are perfect codes of $G$ and $G/N$, respectively, then $H$ is a perfect code of $G$.
\end{enumerate}}
\end{manualtheorem}

\begin{proof}
(a) Suppose that $H$ is a perfect code of $G$. By Lemma \ref{Ct}, there exists a Cayley transversal $T$ of $H$ in $G$. By the definition of a Cayley transversal, we then have $T^{-1}=T$, $G=TH$ and $xH\neq yH$ for any pair of distinct elements $x,y\in T$. Write $T/N = \{xN: x\in T\}$. Then $(T/N)^{-1}=T/N$ and $G/N=(T/N)(H/N)$. Moreover, since $N$ is contained in $H$, we have $(xN)H/N\neq (yN)H/N$ for any pair of distinct elements $xN,yN \in T/N$. Therefore, $T/N$ is a Cayley transversal of $H/N$ in $G/N$. By Lemma \ref{Ct}, $H/N$ is a perfect code of $G/N$.

(b) Suppose to the contrary that $H$ is not a perfect code of $G$. By Theorem \ref{basic}, there exists a $2$-element $x\in G\setminus H$ such that $x^{2}\in H$, $|H:H\cap xHx^{-1}|$ is odd, and $HxH$ contains no involution. Since $N$ is a normal subgroup of $G$ contained in $H$, $xN$ is a $2$-element in $(G/N)\setminus(H/N)$ and $(xN)^{2}\in H/N$. Therefore, $xN(H/N)=x^{-1}N(H/N)$ and it follows that the double coset $(H/N)(xN)(H/N)$ is inverse-closed.
Moreover, $(H/N)(xN)(H/N)$ contains an odd number of left cosets of $H/N$ in $G/N$ as this number is equal to $|H/N:(H\cap xHx^{-1})/N|=|H:H\cap xHx^{-1}|$. Since $H/N$ is a perfect code of $G/N$, by Theorem \ref{basic}, $(H/N)(xN)(H/N)$ contains an involution. Therefore, $(aNxNbN)^{2}=N$ for some $aN,bN\in H/N$. That is, $(axb)^{2}N=N$ and hence $(axb)^{2}\in N$. Since $N$ is a perfect code and a normal subgroup of $G$, by Lemma \ref{HXZ}, there exists $c\in N$ such that $(axbc)^{2}=1$. Note that $a,bc\in H$. However, $HxH$ contains no involution, a contradiction.
\end{proof}

In \cite{ZZ21}, the proof of Theorem 4.2 relies on Corollaries 3.3 and 3.4. So this proof needs to be revised, but fortunately the theorem itself remains true under the same hypothesis. We re-present this theorem below for convenience, and we give a proof of it using Lemma \ref{conjugate} and the corrected Theorem \ref{basic} and Corollary \ref{equivalent2}.

\begin{manualtheorem}{4.2}
\label{ns}
{\em Let $G$ be a group and $H$ a subgroup of $G$. If there exists a Sylow $2$-subgroup of $H$ which is a perfect code of $G$, then $H$ is a perfect code of $G$.}
\end{manualtheorem}

\begin{proof}
Let $g$ be an arbitrary element of $G\setminus H$ such that the double coset $HgH$ is inverse-closed and contains $m$ left cosets of $H$ in $G$ where $m$ is odd. Then $m=|H:H\cap gHg^{-1}|$. In the following we will prove that $HgH$ contains an involution. Once this is achieved, it then follows from Theorem \ref{basic} that $H$ is a perfect code of $G$.

Since $Hg^{-1}H=HgH$, we have $g^{-1}=h_{1}gh_2$ for some $h_1,h_2\in H$. Since $g\notin G\setminus H$ by our assumption, we have $gh_1\notin G\setminus H$. Since $(gh_1)^2=h_{2}^{-1}h_1\in H$, it follows that $gh_{1}$ is of even order. Let $s$ be the largest odd divisor of the order of $gh_{1}$. Set $x=(gh_1)^{s}$. Then $x$ is a $2$-element in $G\setminus H$ with $x^{2}\in H$. Since $x=gh_1(gh_1)^{s-1}$ and $h_1(gh_1)^{s-1}\in H$, we have $HxH=HgH$ and $xHx^{-1}=gHg^{-1}$. Set $L=H\cap xHx^{-1}$. Then $|H:L|=|H:H\cap gHg^{-1}|=m$. Since $xLx^{-1}=xHx^{-1}\cap x^2Hx^{-2}=xHx^{-1}\cap H =L$, $L$ is normal and is of index $2$ in $\langle x,L\rangle$. Let $P$ be a Sylow $2$-subgroup of $\langle x,L\rangle$ such that $x\in P$. Set $Q=P\cap L$. Then $x^{2}\in Q$ and $Q$ is of index $2$ in $P$. It follows that $Q$ is normal in $P$ and is a Sylow $2$-subgroup of $L$. Since $|H:L|=m$ is odd, $Q$ is also a Sylow $2$-subgroup of $H$.
By our assumption, $H$ has a Sylow $2$-subgroup which is a perfect code of $G$. On the other hand, by Sylow's Theorem, any two Sylow $2$-subgroups of $H$ are conjugate in $H$. Thus, by Lemma \ref{conjugate}, $Q$ is a perfect code of $G$. Since $Q$ is normal in $P$ and $x\in P$, we have $x\in N_{G}(Q)$. Hence, by Corollary \ref{equivalent2}, there exists $b\in Q$ such that $xb$ is an involution. Since $Q$ is a subgroup of $H$ and $HxH=HgH$, the involution $xb$ is contained in $HgH$, as required.
\end{proof}

Finally, although \cite[Theorem 4.3]{ZZ21} remains true, its proof needs to be revised due to the reliance on \cite[Corollary 3.3]{ZZ21}. We now give a proof of this result using Theorem \ref{ns} and the corrected Theorem \ref{basic} and Corollary \ref{equivalent2}. We repeat this result here for convenience.

\begin{manualtheorem}{4.3}
\label{QK}
{\em Let $G$ be a group. Let $Q$ be a $2$-subgroup of $G$ and $K$ a $2'$-subgroup of $G$. Suppose that all Sylow $2$-subgroups of $N_{G}(Q)$ are contained in $N_{G}(K)$. Then $QK$ is a perfect code of $G$ if and only if $Q$ is a perfect code of $G$.

In particular, if $Q$ is a $2$-subgroup of $G$ and $K$ is a normal $2'$-subgroup of $G$, then $QK$ is a perfect code of $G$ if and only if $Q$ is a perfect code of $G$.}
\end{manualtheorem}

\begin{proof}
Since all Sylow $2$-subgroups of $N_{G}(Q)$ are contained in $N_{G}(K)$, $Q$ is contained in $N_{G}(K)$. Hence $H=QK$ is a subgroup of $G$ and $K$ is normal in $H$.

Since the order of $Q$ is a power of $2$ but the order of $K$ is odd, we have $Q\cap K=\{1\}$. Therefore, $Q$ is a Sylow $2$-subgroup of $H$ and $K$ is a Hall $2'$-subgroup of $H$. Thus, by Theorem \ref{ns}, if $Q$ is a perfect code of $G$, then $H$ is a perfect code of $G$. This proves the sufficiency.

We now prove the necessity. Assume that $H$ is a perfect code of $G$. Consider an arbitrary $2$-element $x \in N_{G}(Q)$ with $x^{2}\in Q$. Since $x^{2}\in H$, we have $(HxH)^{-1}=HxH$. Since all Sylow $2$-subgroups of $N_{G}(Q)$ are contained in $N_{G}(K)$, we have $x\in N_{G}(K)$. Therefore, $HxH=QKxH=xQKH=xH$.
By Theorem \ref{basic}, $HxH$ contains an involution. Since $HxH=xQK$,
we have $(xab)^{2}=1$ for some $a\in Q$ and $b\in K$.
Noting that $xa\in N_{G}(K)$, we obtain $(xa)^{2}=xab^{-1}a^{-1}x^{-1}b^{-1}\in K$. On the other hand, since $x\in N_{G}(Q)$, we get $(xa)^{2}=x^{2}x^{-1}axa\in Q$. Recalling $Q\cap K=\{1\}$, it follows that $(xa)^{2}=1$. Thus, by Corollary \ref{equivalent2}, $Q$ is a perfect code of $G$.
\end{proof}

\smallskip

\noindent {\textbf{Acknowledgements}}~~
\medskip

We would like to thank Dr. Kai Yuan for bringing an error in the proof of \cite[Theorem 3.1]{ZZ21} to our attention. We are grateful to the anonymous referees for their careful reading and valuable suggestions on this paper.  The first author was supported by the Basic Research and Frontier
Exploration Project of Chongqing (cstc2018jcyjAX0010), the Science and Technology Research
Program of Chongqing Municipal Education Commission (KJQN201800512) and the National Natural
Science Foundation of China (No. 11671276). The second author was supported by the Research
Grant Support Scheme of The University of Melbourne.





\end{document}